
\mag=\magstep1
\documentstyle{amsppt}

\topmatter
\title On Fujita's freeness conjecture for 3-folds and 4-folds
\endtitle
\author Yujiro Kawamata
\endauthor

\rightheadtext{Fujita's freeness conjecture}

\address Department of Mathematical Sciences, University of Tokyo, Komaba,
Meguro, Tokyo, 153, Japan \endaddress
\email kawamata\@tansei.cc.u-tokyo.ac.jp\endemail


\abstract
We shall prove a conjecture of T. Fujita on the freeness of the
adjoint linear systems in some cases:
Let $X$ be a smooth projective variety of dimension $n$
and $H$ an ample divisor.  Assume that $n = 3$ or $4$.
Then $\vert K_X + mH \vert$ is free if $m \ge n+1$.
Moreover, we obtain more precise result in the case $n = 3$.
\endabstract

\endtopmatter

\document

\head Introduction
\endhead

T. Fujita raised the following:

\proclaim{Conjecture 1}
Let $X$ be a smooth projective variety of dimension $n$ and
$H$ an ample divisor.
Then $\vert K_X + mH \vert$ is free if $m \ge n+1$.  Moreover, if
$(H^n) \ge 2$, then $\vert K_X + nH \vert$ is also free.
\hfill $\square$
\endproclaim

In the case $n = 3$, Ein and Lazarsfeld [EL1]
gave an affirmative answer to the first part of
Conjecture 1 and Fujita [F] the second part.
A stronger version of Fujita's freeness conjecture is the following:

\proclaim{Conjecture 2}
Let $X$ be a normal projective variety of dimension $n$, $x_0 \in X$
a smooth point, and $L$ an ample Cartier divisor.
Assume that there exist positive numbers $\sigma_p$ for
$p = 1,2, \ldots, n$ which satisfy the following conditions:

(1) $\root{p}\of{(L^p \cdot W)} \ge \sigma_p$ for any subvariety $W$ of
dimension $p$ which contains $x_0$,

(2) $\sigma_p \ge n$ for all $p$ and $\sigma_n > n$.

Then $\vert K_X + L \vert$ is free at $x_0$.
\hfill $\square$
\endproclaim

In the case $n = 3$, [F] proved that, if
$\sigma_1 \ge 3, \sigma_2 \ge \sqrt 7$ and $\sigma_3 \ge \root{3}\of{51}$,
then $\vert K_X + L \vert$ is free at $x_0$.
In an arbitrary dimension,
Angehrn and Siu [AS] proved a weaker result that, if
$\sigma_p > \frac12 n(n+1)$ for any $p$,
then $\vert K_X + L \vert$ is free at $x_0$.
In particular, $\vert K_X + mH \vert$ is free if $m \ge \frac 12 n(n+1) +1$
in the situation of Conjecture 1.
Parallel arguments are possible in differential geometry and
algebraic geometry;
this proof is translated to algebraic geometry by Koll\'ar.
There is also a paper by Tsuji [T].
A paper of Smith [Sm] suggests that, even if $X$ is singular,
we should have the spannedness of the reflexive sheaf $\Cal O_X(K_X + L)$
under certain conditions.

We shall prove the following results in this paper:

(1) (Theorem 3.1): In the case $n = 3$, Conjecture 2 is true.

(2) (Theorem 4.1): In the case $n = 4$, if $\sigma_p \ge 5$ for all $p$, then
$\vert K_X + L \vert$ is free at $x_0$.
In particular, the first part of Conjecture 1 is true.

In \S\S 1 and 2,
we shall explian the general strategy toward the freeness results.
This is an application of the vanishing theorem of [K1] and [V], and
has the origin in the proof of the base point free theorem ([K2], [Sh1]).
The adjoint linear system appears naturally in the course of the proof.
We note that there exists no universal bound of $m_0$ as a function of $n$
such that $\vert mH \vert$ is free if $m \ge m_0$.
This justifies that
we ask the effective freeness for the adjoint linear system
$\vert K_X + mH \vert$.
\S\S 3 and 4 are devoted to the proof of the main results.

The author would like to thank Professors L. Ein and R. Lazarsfeld
for showing the manuscript [EL2] to the author at AMS Summer Institute at
Santa Cruz.
In the first version of this paper, Theorem 2.2 was weaker so that
we had to assume
$\sigma_3 > \root{3}\of{3(\frac2{3-\sqrt2})^3 + 24}
= \root{3}\of{30.0183 \cdots} = 3.107865 \cdots$
in Theorem 3.1.  After that,
the author received two letters, one from S. Helmke and
the other from L. Ein and R. Lazarsfeld, and both contained the optimal result
stated as in Theorem 3.1 by using the idea of T. Fujita [F].
This paper follows the argument of the former.
The author would like to express his gratitude to S. Helmke
for allowing the author to reproduce his result.

\head 1. Minimal center of log canonical singularities
\endhead

We recall the standard notation (cf. [KMM]).
Let $X$ be a normal variety of dimension $n$.
A $\Bbb Q$-{\it divisor}
is an element of $Z_{n-1}(X) \otimes \Bbb Q$, i.e., a finite formal sum
$D = \sum_{j=1}^n d_jD_j$ of prime divisors $D_j$ with coefficients
$d_j \in \Bbb Q$.  We usually require implicitly
that the $D_j$ are distinct.
$D$ is said to be {\it effective} if $d_j \ge 0$ for all $j$.
The {\it round up} of $D$ is defined by
$\ulcorner D \urcorner = \sum_j \ulcorner d_j \urcorner D_j$.
Two $\Bbb Q$-divisors $D_1, D_2$ are said to be $\Bbb Q$-linearly equivalent,
and we write $D_1 \sim_{\Bbb Q} D_2$, if there exists a positive integer
$m$ and a non-zero rational function $h$ such that
$m(D_1 - D_2) = \text{div}(h)$.  By abuse of notation, we sometimes write
$D_1 = D_2$ instead of $D_1 \sim_{\Bbb Q} D_2$
when the canonical divisors are involved
(e.g., the first paragraph of Definition 1.2).

$D$ is called a $\Bbb Q$-{\it Cartier divisor} if
it is in the image of the natural injective homomorphism
$\text{Div}(X) \otimes \Bbb Q \to Z_{n-1}(X) \otimes \Bbb Q$.
If $D$ is an effective $\Bbb Q$-Cartier divisor and $x_0 \in X$ is a point,
then the order $\text{ord}_{x_0} D \in \Bbb Q$ is defined by linearity.
If $X$ is complete and $D$ is $\Bbb Q$-Cartier,
we can define the intersection number $(D^s \cdot S) \in \Bbb Q$
for any subvariety $S$ of dimension $s$ on $X$.
$D$ is said to be {\it nef} if $(D \cdot C) \ge 0$ for any curve $C$.
In this case, $D$ is called {\it big} if $(D^n) > 0$.

Let $\mu: Y \to X$ be a birational morphism of normal varieties.
The {\it exceptional locus} $\text{Exc}(\mu)$ of $\mu$
is the smallest closed subset of
$Y$ such that $\mu \vert_{Y \setminus \text{Exc}(\mu)}$ is an isomorphism.
If $D$ is a $\Bbb Q$-divisor
(resp. a $\Bbb Q$-Cartier divisor) on $X$,
we can define the {\it strict transform} (resp. {\it total transform})
$\mu_*^{-1}D$ (resp. $\mu^*D$).
For a pair $(X, D)$ of a variety and a $\Bbb Q$-divisor, an
{\it embedded resolution} or a {\it log resolution} is a
proper birational morphism $\mu: Y \to X$ from a smooth variety $Y$ such that
the union of the support of $\mu_*^{-1}D$ and $\text{Exc}(\mu)$ is
a normal crossing divisor.

Most of the results of this paper are the
applications of the following vanishing theorem:

\proclaim{Theorem 1.1}
Let $X$ be a smooth projective variety and $D$ a $\Bbb Q$-divisor.
Assume that $D$ is nef and big, and that the support of the difference
$\ulcorner D \urcorner - D$ is a normal crossing divisor.
Then $H^p(X, K_X + \ulcorner D \urcorner) = 0$ for $p > 0$.
\hfill $\square$
\endproclaim

\definition{Definition 1.2}
Let $X$ be a normal variety and $D = \sum_id_iD_i$ an effective
$\Bbb Q$-divisor such that $K_X + D$ is $\Bbb Q$-Cartier.
If $\mu: Y \to X$ is an embedded resolution of the pair $(X, D)$,
then we can write
$$
K_Y + \mu^{-1}_*D = \mu^*(K_X + D) + F
$$
with $F = \sum_j e_jE_j$ for the exceptional divisors $E_j$.
We call $F$ the {\it discrepancy} and $e_j \in \Bbb Q$ the
{\it discrepancy coefficient} for $E_j$.  We regard $- d_i$ as the
discrepancy coefficient for $D_i$.

The pair $(X, D)$ is said to have only {\it log canonical singularities (LC)}
(resp. {\it kawamata log terminal singularities (KLT)})
if $d_i \le 1$ (resp. $< 1$) for all $i$ and $e_j \ge -1$ (resp. $>-1$)
for all $j$
for an embedded resolution $\mu: Y \to X$.
One can also say that $(X, D)$ is LC (resp. KLT),
or $K_X + D$ is LC (resp. KLT), when $(X, D)$ has only
LC (resp. KLT).  If $D = 0$, then $X$ is called LC (resp. KLT) if so is
$(X, D)$.
The pair $(X, D)$ is said to be LC (resp. KLT) at a point $x_0 \in X$
if $(U, D \vert_U)$ is LC (resp. KLT) for some neighborhood $U$ of $x_0$.
\enddefinition

\definition{Definition 1.3}
A subvariety $W$ of $X$ is said to be a
{\it center of log canonical singularities} for the pair $(X, D)$,
if there is a birational morphism from a normal variety
$\mu: Y \to X$ and a prime divisor $E$ on $Y$
with the discrepancy coefficient $e \le - 1$
such that $\mu(E) = W$.
For example, if $E = \mu_*^{-1}D_i$ for some $i$,
then we have $e = - d_i$, so $D_i$ is a center of log canonical
singularities if and only if $d_i \ge 1$.
For another such $\mu': Y' \to X$, if the strict transform $E'$ of $E$ exists
on $Y'$, then we have the same discrepancy coefficient for $E'$.
The divisor $E'$ is considered to be equivalent to $E$, and
the equivalence class of these prime divisors is called
a {\it place of log canonical singularities} for $(X, D)$.

The set of all the centers (resp. places) of log canonical singularities
is denoted by $CLC(X, D)$ (resp. $PLC(X, D)$).
Thus there is a natural surjective map
$PLC(X, D) \to CLC(X, D)$, which is
not necessarily injective.
If $(X, D)$ is LC, then $CLC(X, D)$ is a finite set.
The union of all the subvarieties in $CLC(X, D)$ is denoted by $LLC(X, D)$
and called the {\it locus of log canonical singularities} for $(X, D)$.
$LLC(X, D)$ is a closed subset of $X$, and is empty if and only if
$(X, D)$ is KLT.
For a point $x_0 \in X$, we define $CLC(X, x_0, D)
= \{W \in CLC(X, D) ; x_0 \in W\}$.
\enddefinition

\proclaim{Theorem 1.4} ({\it Connectedness Lemma},[Sh2], [Ko]).
Let $f: X \to Z$ be a proper surjective morphism of normal varieties
with connected fibers, and
$D = \sum _i d_iD_i$
a $\Bbb Q$-divisor on $X$ such that $K_X + D$ is $\Bbb Q$-Cartier.
Assume the following conditions:

\roster
\item if $d_i < 0$, then $\text{codim}(f(D_i)) \ge 2$,

\item $- (K_X + D)$ is $f$-nef and $f$-big.
\endroster

Then $LLC(X, D) \cap f^{-1}(z)$ is connected
for any point $z \in Z$.
\endproclaim

We include the proof for the convenience of the reader.

\demo{Proof}
Let $\mu: Y \to X$ be an embedded resolution of the pair $(X, D)$, and
$K_Y + D_Y = \mu^*(K_X + D)$.  By definition,
we have $LLC(X, D) = \mu(LLC(Y, D_Y))$.
So we may assume that $X$ is smooth and $\text{Supp}(D)$ is normal
crossing.
If we write $D = S + D'$ with
$S = \sum _{d_i\ge1} d_iD_i$ and
$D' = \sum _{d_i<1} d_iD_i$, then
$LLC(X, D) = \text{Supp}(S)$.
By the condition (2),
we have $R^1f_*\Cal O_X(\ulcorner - D \urcorner) = 0$.
{}From an exact sequence
$$
0 \to \Cal O_X(\ulcorner - D \urcorner)
\to \Cal O_X(\ulcorner - D' \urcorner)
\to \Cal O_{\llcorner S \lrcorner}(\ulcorner - D' \urcorner) \to 0
$$
we deduce that the natural homomorphism
$f_*\Cal O_X(\ulcorner - D' \urcorner)
\to f_*\Cal O_{\llcorner S \lrcorner}(\ulcorner - D' \urcorner)$
is surjective.
Since $\ulcorner - D' \urcorner \ge 0$, there is a natural homomorphism
$\Cal O_Z \to f_*\Cal O_X(\ulcorner - D' \urcorner)$,
which is an isomorphism by the condition (1).
Therefore, we obtain our assertion.
\hfill $\square$ \enddemo

\proclaim{Proposition 1.5}
Let $X$ be a normal variety and $D$ an effective
$\Bbb Q$-Cartier divisor such that $K_X + D$ is $\Bbb Q$-Cartier.
Assume that $X$ is KLT and $(X, D)$ is LC.
If $W_1, W_2 \in CLC(X, D)$ and $W$ an irreducible component of
$W_1 \cap W_2$, then $W \in CLC(X, D)$.
In particular, if $(X, D)$ is not KLT at a point $x_0 \in X$,
then there exists the unique minimal element of $CLC(X, x_0, D)$.
\endproclaim

\demo{Proof}
Since the assertion is local, we may assume that $X$ is affine.
Let $D_i$ ($i=1,2$) be a general member among effective Cartier divisors
which contain $W_i$.
Let $\mu: Y \to X$ be an embedded resolution
of the pair $(X, D + D_1 + D_2)$.
We choose $\mu$ so that there are
divisors $E_i$ above the $W_i$ with the discrepancy coefficients $-1$.
Let $e_i$ (resp. $e_i'$) be the coefficients of $E_i$ in $\mu^*D$
(resp. $\mu^*D_i$), and $a_i$ the positive numbers such that
$e_i = a_ie'_i$.  Then we have
$LLC(X, (1-\epsilon)D + a_1\epsilon D_1 + a_2\epsilon D_2) = W_1 \cup W_2$
for $0 < \epsilon \ll 1$.
By the connectedness lemma,
there exist divisors $F_i(\epsilon)$ on $Y$ for $i = 1,2$ such that
$F_i(\epsilon) \in PLC(X, (1-\epsilon)D + a_1\epsilon D_1 + a_2\epsilon D_2)$,
$\mu(F_i(\epsilon)) \subset W_i$, and $F_1(\epsilon) \cap F_2(\epsilon) \ne
\emptyset$.
Since there are only a finite number of the exceptional divisors
for $\mu$, we have the common divisors $F_i = F_i(\epsilon)$
for a convergent sequence
$\epsilon \to 0$.  Then $\mu(F_1 \cap F_2) \in CLC(X, D)$.
\hfill $\square$ \enddemo

We shall control the singularities of the minimal center of log canonical
singularities in the following theorem which
are also obtained independently by Ein and Lazarsfeld ([EL2]).
This is a variant of the connectedness lemma.

\proclaim{Theorem 1.6}
Let $X$ be a normal vairiety,
$x_0 \in X$ a KLT point,
$D$ an effective $\Bbb Q$-Cartier divisor such that
$(X, D)$ is LC at $x_0$, $W$ a minimal element of
$CLC(X, x_0, D)$,
and $E$ a place of log canonical singularities for
$(X, D)$ on a blow-up $Y$ of $X$ which lies above $W$.
Then $W$ is normal at $x_0$, and
the projection $E \to W$ has connected fibers in a neighborhood of $x_0$.
\endproclaim

\demo{Proof}
Since the assertions are local, we may assume that $X$ is affine.
We may assume that $LLC(X, D) = W$ as in the proof of Proposition 1.5.
We may also assume that $PLC(X, D)$ consists of one element.
Let $\mu: Y \to X$ be an embedded resolution of the pair $(X, D)$.
We write
$$
K_Y + E + F = \mu^*(K_X + D)
$$
where $E$ is a prime divisor such that $\mu(E) = W$, and
$F$ is a divisor whose coefficients are smaller than 1.
By Theorem 1.1
$$
H^1(Y, - E + \ulcorner - F \urcorner) = 0
$$
and we obtain a surjection
$$
H^0(Y, \ulcorner - F \urcorner) \to H^0(E, \ulcorner - F \urcorner \vert_E)
$$
Since $\ulcorner - F \urcorner$ is effective and exceptional,
we have
$$
H^0(X, \Cal O_X) @>{\sim}>> H^0(Y, \ulcorner - F \urcorner)
$$
Therefore, the natural injective homomorphism
$$
H^0(W, \Cal O_W) \to H^0(E, \ulcorner - F \urcorner)
$$
is surjective.  Thus $W$ is normal, and $\mu: E \to W$ has connected fibers.
\hfill $\square$
\enddemo

\proclaim{Proposition 1.7}
Let $x_0 \in X$, $D$ and $W$ be as in Theorem 1.6.
Assume that $W$ is a prime divisor.
Then there exists an effective $\Bbb Q$-divisor $D_W$ on $W$
such that
$K_X + D \vert_W = K_W + D_W$ and
the pair $(W, D_W)$ is KLT at $x_0$.
\endproclaim

\demo{Proof}
By using the residue map $\omega_X(W) \to \omega_W$,
we can define naturally the effective $\Bbb Q$-divisor
$D_W$ on $W$ (cf. [KMM]).
We have
$$
\mu_E^*(K_W + D_W) = (K_Y + E + F) \vert_E = K_E + F \vert_E
$$
where $\mu_E = \mu \vert_E: E \to W$,
and the pair $(W, D_W)$ is KLT, since $\mu_E$ is birational.
\hfill $\square$
\enddemo

\definition{Question 1.8}
Let $x_0 \in X$, $D$ and $W$ be as in Theorem 1.6.
Does there exist an effective $\Bbb Q$-divisor $D_W$ on $W$
such that
$K_X + D \vert_W = K_W + D_W$
and the pair $(W, D_W)$ is KLT at $x_0$?
\enddefinition

We shall give an affirmative answer to Question 1.8 in the case
$\text{codim }_X W = 2$ in [K3].
We also have a positive evidence in the case $\text{dim }W = 2$:

\proclaim{Theorem 1.9}
Let $x_0 \in X$, $D$ and $W$ be as in Theorem 1.6.
Assume that $\text{dim }W = 2$.
Then $W$ has at most a rational singularity at $x_0$.
Moreover, if $W$ is singular at $x_0$,
and if $D'$ is an effective $\Bbb Q$-Cartier divisor
on $X$ such that $\text{ord}_{x_0} D' \vert_W \ge 1$, then
$\{x_0\} \in CLC(X, x_0, D + D')$.
\endproclaim

\demo{Proof}
We use the notation of Theorem 1.6.
As in the proof of Proposition 1.5,
we may assume that none of the coefficients of $F$ are integers
and the projection $\mu: E \to W$ is factorized as
$E @>{\alpha}>> \tilde W @>{\beta}>> W$ with $\tilde W$ smooth and
$\beta$ birational.
By the vanishing theorem
$$
H^p(Y, - E + \ulcorner - F \urcorner) =
H^p(Y, \ulcorner - F \urcorner) = 0 \text{ for } p > 0
$$
where the latter is obtained by replacing $D$ by
$(1-\epsilon)D$ for a small $\epsilon > 0$.
Then
$$
H^p(E, \ulcorner - F \vert_E \urcorner) = 0 \text{ for } p > 0
$$
Similarly, we have $R^p\alpha_*\Cal O_E(\ulcorner - F \vert_E \urcorner) = 0$
for $p > 0$.
Therefore, if we set
$\Cal F = \alpha_*\Cal O_E(\ulcorner - F \vert_E \urcorner)$,
then $R^p\beta_*\Cal F = 0$ for $p > 0$.
We know that $\ulcorner - F \vert_E \urcorner$
is effective and $\mu_*\Cal O_E(\ulcorner - F \vert_E \urcorner) = \Cal O_W$.
Therefore, we have
an injective homomorphism $\Cal O_{\tilde W} \to \Cal F$, and
the support of its cokernel is contained in
$\text{Exc}(\beta) = \bigcup_i G_i$.

Let $\Cal F^{**} = \Cal O_{\tilde W}(G)$ be the double dual of
$\Cal F$, where $G = \sum_i g_iG_i$ for some non-negative integers $g_i$.
Let $g = \sum_i g_i$.  We define a sequence of effective divisors
$G^j = \sum_i g_i^jG_i$ for $j = 0, 1, \ldots, g$
such that $\sum_i g_i^j = g - j$ inductively as follows.
First, set $G^0 = G$.
Assuming that $G^{j_0}$ is defined for a $j_0 < g$, pick an $i_{j_0}$
such that $(G^{j_0} \cdot G_{i_{j_0}}) < 0$, and
set $G^{j_0+1} = G^{j_0} - G_{i_{j_0}}$.
Let $\Cal F^j = \Cal F \cap \Cal O_{\tilde W}(G^j)$.
Then we have injective homomorphisms
$\Cal F^j/\Cal F^{j+1} \to \Cal O_{\tilde W}(G^j)/\Cal O_{\tilde W}(G^{j+1})
\cong \Cal O_{G_{i_j}}(G^j)$.
Since $H^0(G_{i_j}, G^j) = 0$, we have
$H^1(\Cal F^{j+1}) = 0$ if $H^1(\Cal F^j) = 0$.  Therefore, we have
$H^1(\tilde W, \Cal O_{\tilde W}) = 0$.

In order to prove the second part, we shall prove that
$g_{i_0} = 0$ for some $i_0$.
Assume the contrary that $g_i > 0$ for all $i$.
By the above argument, we have
$H^1(\Cal F^j/\Cal F^{j+1}) = 0$, hence
$H^1(\Cal O_{\tilde W}(G^j)/\Cal O_{\tilde W}(G^{j+1})) = 0$
for all $j$.
Then $H^1(G, G \vert _G) = 0$.
By the Serre duality, we have $H^0(G, K_{\tilde W} \vert_G) = 0$.
But since $W$ has a rational singularity at $x_0$,
it is a contradiction, hence $g_{i_0} = 0$ for some $i_0$.

Let $F \vert_E = \sum_{\ell} f_{\ell}F_{\ell}$ and
$\alpha^*G_{i_0} = \sum_{\ell} k_{\ell}F_{\ell}$.
Since $g_{i_0} = 0$,
there exists an $\ell_0$ such that
$\alpha(F_{\ell_0}) = G_{i_0}$ and
$\ulcorner - f_{\ell_0} \urcorner < k_{\ell_0}$,
hence $-f_{\ell_0} \le k_{\ell_0} - 1$.
We have $K_Y + E + F + \mu^*D' = \mu^*(K_X + D + D')$.
Since $\text{ord}_{x_0} D' \vert_W \ge 1$, the coefficient of $F + \mu^*D'$
on $F_{\ell_0}$ is at least $f_{\ell_0} + k_{\ell_0} \ge 1$,
and $\{x_0\} \in CLC(X, x_0, D + D')$.
\hfill $\square$
\enddemo

We shall replace the minimal center of log canonical singularities
by a smaller subvariety by using the following theorem due to [EL2].
This is another evidence of Question 1.8 on the adjunciton and the
inverse adjunction.

\proclaim{Theorem 1.10}
Let $x_0 \in X$, $D$ and $W$ be as in Theorem 1.6.
Let $D_1$ and $D_2$ be effective $\Bbb Q$-Cartier divisors on $X$
whose supports do not contain $W$ and which induce the same
$\Bbb Q$-Cartier divisor on $W$.
Assume that $(X, D + D_1)$ is LC at $x_0$ and
there exists an element of $CLC(X, x_0, D + D_1)$ which is
properly contained in $W$.
Then the similar statement holds for the pair $(X, D + D_2)$.
\endproclaim

\demo{Proof}
Since the assertion is local, we may assume that $X$ is affine.
First, we assume that $W$ is the only element of $CLC(X, D)$ and
there is only one place of log canonical singularities above it.
Let $\mu: Y \to X$ be an embedded resolution of the pair
$(X, D + D_1 + D_2)$, and $E$ the divisor which represents
the place of log canonical singularities for $(X, D)$.
We write $K_Y + E + \sum_i f_iF_i = \mu^*(K_X + D)$ and
$K_Y + E + \sum_i f_{i,\alpha}F_i = \mu^*(K_X + D + D_{\alpha})$
for $\alpha = 1,2$.
Let $D_W$ be the $\Bbb Q$-Cartier divisors on $W$ which is induced by the
$D_{\alpha}$, and $\mu^*D_W = \sum_i g_iG_i$ for $G_i = F_i \cap E$.
Then we have $f_{i,1} = f_{i,2} = f_i + g_i$ if $G_i \ne \emptyset$.
By the connectedness lemma,
there exists $i_0$ such that $G_{i_0} \ne \emptyset$ and $f_{i_0,1} \ge 1$.
Since the support of $D_1$ does not contain $W$, we have
$\mu(G_{i_0}) \subsetneq W$.
Hence there exists an element of $CLC(X, x_0, D + D_2)$ which is
properly contained in $W$.

Let
$$
c = \sup \{t \in \Bbb Q ; K_X + D + tD_2 \text{ is LC at }x_0 \}
$$
Then the assumption of the theorem is satisfied by the pair
$(X, D + cD_2)$.
By the preceding argument, there exists an element of
$CLC(X, x_0, D + cD_1)$ which is
properly contained in $W$.
Since $(X, D + D_1)$ is LC, we have
$c \ge 1$, and $(X, D + D_2)$ is LC at $x_0$.

Now we consider the general case.
Let us take a general Cartier divisor $D'$ which contains $W$, and
a positive number $a$ such that
$K_X + (1-\epsilon)D + a\epsilon D'
= K_X + D + \epsilon(-D + aD')$ is LC with the only one
place of canonical singularities for any $0 < \epsilon \ll 1$.
There exists a function $b(\epsilon)$ such that
the assumption of the theorem holds for the pair
$(X, (1-\epsilon)D + a\epsilon D' + (1+b(\epsilon))D_1)$.
Being LC is a closed condition for the coefficients of
divisors, so $b(\epsilon)$ is a well defined continuous convex function
such that $\lim_{\epsilon \to 0} b(\epsilon) = 0$.
By the first part of the proof, the conclusion of the theorem
holds for $(X, (1-\epsilon)D + a\epsilon D' + (1+b(\epsilon))D_2)$.
Looking at the discrepancies, we conclude that $c = 1$.
Then the minimal element of
$CLC(X, x_0, D + D_2)$ should be smaller than $W$,
since the support of $D_2$ does not contain $W$.
\hfill $\square$
\enddemo

\head 2. General method
\endhead

We shall consider the conditions for the existence of a member of the
given linear system which has prescribed order at a given point.

\proclaim{Proposition 2.1}
Let $X$ be a normal and complete variety of dimension $n$,
$L$ a nef and big Cartier divisor,
$x_0 \in X$ a point, and $t$ a rational number such that $t > 1$.
Then there exists an effective $\Bbb Q$-Cartier divisor $D$ such that
$D \sim_{\Bbb Q} tL$ and
$$
\text{ord}_{x_0} D \ge \root{n}\of{\frac{(L^n)}{\text{mult}_{x_0} X}}
$$
\endproclaim

\demo{Proof}
We shall prove that
there exists a positive integer $m$ with $mt \in \Bbb N$ and a member
$D_m \in \vert mtL \vert$ such that
$$
\text{ord}_{x_0} D_m \ge m\root{n}\of{\frac{(L^n)}{\text{mult}_{x_0} X}}
$$
Since
$$
\text{length }\Cal O_{X,x_0}/\frak m_{x_0}^d
= \frac{d^n}{n!}\text{mult}_{x_0} X + \text{lower terms in }d
$$
it is enough to prove that
$$
h^0(X, mL) = \frac{m^n(L^n)}{n!} + \text{lower terms in }m
$$
If we replace $X$ by its desingularization, we may assume that $X$ is smooth
and projective.
Let $H$ be a very ample Cartier divisor such that $H - K_X$ is
ample.  Then we have
$H^p(X, mL + H) = 0$ for any $p>0$ and $m > 0$.
Hence
$$
h^0(X, mL + H) = \chi(X, mL + H) =
\frac{m^n(L^n)}{n!} + \text{ lower terms in }m
$$
If $Y$ is a general member of the linear system $\vert H \vert$, then
$$
0 \to \Cal O_X(mL) \to \Cal O_X(mL + H) \to \Cal O_Y(mL + H) \to 0
$$
Since $\text{dim }Y = n - 1$, we obtain the result.
\hfill $\square$
\enddemo

The following theorem is due to S. Helmke after the idea of Fujita [F].

\proclaim{Theorem 2.2}
Let $X$ be a normal projective variety of dimension $n$,
$x_0 \in X$ a smooth point, $L$ an ample Cartier divisor, $W$ a prime
divisor with $\text{ord}_{x_0} W = d \ge 1$,
and $e, k$ positive rational numbers
such that $de \le 1$ and $k^n < (L^n)$.
Assume that for any effective $\Bbb Q$-divisor $D$,
if $D \sim_{\Bbb Q} L$ and $\text{ord}_{x_0} D \ge k$,
then it follows that $D \ge ekW$.
Then there exists a real number $\lambda$ with $0 \le\lambda \le 1 - de$
which satisfies the following condition:
if $k'$ is a positive rational number such that $k' > k$ and
$$
(\lambda k)^n + (\frac{1 - de - \lambda}{1 - \lambda})^{n-1}
\{(k' + \frac{\lambda dek}{1 - de - \lambda})^n
- (\lambda k + \frac{\lambda dek}{1 - de - \lambda})^n\} < (L^n)
$$
then there exist an effective $\Bbb Q$-divisor $D$
such that $D \sim_{\Bbb Q} L$ and $\text{ord}_{x_0} D \ge k'$.
(If $\lambda = 1 - de$, then the left hand side of the above inequality
should be taken as a limit.)
\endproclaim

\demo{Proof}
Let us define a function $\phi(q)$ for $q \in \Bbb Q_{\ge 0}$ to be the
largest real number such that $D \ge \phi(q)W$ whenever
$D \ge 0, D \sim_{\Bbb Q} L$ and $\text{ord}_{x_0} D = q$.
Then $\phi$ is a convex function.  In fact, if
$\text{ord}_{x_0} D_i = q_i$ and $D_i = (\phi(q_i) + \epsilon_i)W +
\text{other components}$
for $0 < \epsilon_i \ll 1$ and $i = 1,2$,
then $\text{ord}_{x_0} (tD_1 + (1-t)D_2) = tq_1 + (1-t)q_2$ and
$tD_1 + (1-t)D_2 = (t(\phi(q_1) + \epsilon_1) + (1-t)(\phi(q_2)+ \epsilon_2))W
+ \text{other components}$, hence
$\phi(tq_1 + (1-t)q_2) \le t\phi(q_1) + (1-t)\phi(q_2)$.
Since $\phi(k) \ge ek$, there exists a real number
$\lambda$ such that $0 \le \lambda < 1$ and
$\phi(q) \ge \frac{e(q-\lambda k)}{1-\lambda}$ for any $q$.
Since $q \ge \phi(q)d$, we have $\lambda \le 1 - de$.

Let $m$ be a large and sufficiently divisible integer and
$v: H^0(X, mL) \to
\mathbreak
\Cal O_{X,x_0}(mL) \cong \Cal O_{X,x_0}$
the evaluation homomorphism.
We consider subspaces $V_{j} = v^{-1}(\frak m_{x_0}^j)$
of $H^0(X, mL)$ for integers $j$ such that $\lambda km \le j \le k'm$.
First, we have
$$
\text{dim }V_{\ulcorner \lambda km \urcorner} \ge \text{dim }H^0(X, mL)
- \frac{(\lambda km)^n}{n!} + \text{lower terms in }m
$$
Let $D \in \vert mL \vert$ be a member corresponding to
$h \in V_{j}$ for some $j$.
Since we have $D \ge \phi(j/m)mW$,
the number of conditions in order for $h \in V_{j+1}$
is at most the number of homogeneous polynomials of order $j - \phi(j/m)dm$ in
$n$ veriables, i.e.,
$$
\frac{(j - \phi(j/m)dm)^{n-1}}{(n-1)!} + \text{lower terms in }m
$$
Therefore, our assertion follows from
$$
\align
&\frac{(\lambda km)^n}{n!} +
\sum_{j=\ulcorner \lambda km \urcorner}^{k'm-1}
\frac{(j - \frac{de(j-\lambda km)}{1-\lambda})^{n-1}}{(n-1)!}
+ \text{lower terms in }m \\
&=\frac{(\lambda km)^n}{n!} +
\frac{(\frac{1 - de - \lambda}{1 - \lambda})^{n-1}
\{(k' + \frac{\lambda dek}{1 - de - \lambda})^n
- (\lambda k + \frac{\lambda dek}{1 - de - \lambda})^n\}m^n}{n!} \\
&\,\,\,+ \text{lower terms in }m \\
&< \frac{m^n(L^n)}{n!} + \text{lower terms in }m
\endalign
$$
\hfill $\square$
\enddemo

We begin to state our freeness result in the ideal case in which
the minimal center of canonical singularities is an isolated point.

\proclaim{Proposition 2.3}
Let $X$ be a normal projective variety of dimension $n$,
$x_0 \in X$ a Gorenstein KLT point,
and $L$ an ample Cartier divisor.
Assume that there exists an effective $\Bbb Q$-Cartier divisor
$D$ which satisfies the following conditions:

(1) $D \sim_{\Bbb Q} tL$ for a rational number $t < 1$,

(2) $(X, D)$ is LC at $x_0$,

(3) $\{x_0\} \in CLC(X, D)$.

Then $\vert K_X + L \vert$ is free at $x_0$.
\endproclaim

\demo{Proof}
Let $D'$ be a general member of $\vert mL \vert$ for $m \gg 0$
which passes through $x_0$.
If we replace $D$ by $(1 - \epsilon_1)(D + \epsilon_2 D')$
for some $0 < \epsilon_i \ll \frac1m$, then
we may assume that $x_0$ is an isolated point of $LLC(X, D)$.
Let $\mu: Y \to X$ be an embedded resolution of the pair $(X, D)$.
Then
$$
K_Y + E + F = \mu^*(K_X + D)
$$
where $E$ is a reduced divisor such that $\mu(E) = \{x_0\}$, and
$F$ is a divisor of the form $\sum_j f_jF_j$ with $f_j < 1$
if $x_0 \in \mu(F_j)$.
Then
$$
K_Y + (1 - t)\mu^*L
\sim_{\Bbb Q} \mu^*(K_X + L) - E - F.
$$
Thus
$$
H^1(Y, \mu^*(K_X + L) - E + \ulcorner - F \urcorner) = 0
$$
and we obtain a surjection
$$
H^0(Y, \mu^*(K_X + L) + \ulcorner - F \urcorner)
\to H^0(E, \mu^*(K_X + L)) \cong \Bbb C.
$$
Since $\ulcorner - F \urcorner$ is effective and exceptional over
a neighborhood of $x_0$,
$$
H^0(X, K_X + L) \to H^0(E, \mu^*(K_X + L))
$$
is also surjective.
\hfill $\square$
\enddemo

By combining Propositions 2.1 and 2.3, we would obtain Fujita's
freeness conjecture if we would have only the ideal case.

\head 3. Smooth 3-fold
\endhead

\proclaim{Theorem 3.1}
Let $X$ be a normal projective variety of dimension $3$,
$L$ an ample Cartier divisor, and $x_0 \in X$ a smooth point.
Assume that there are positive numbers $\sigma_p$ for $p = 1,2,3$
which satisfy the following conditions:

\roster
\item $\root{p}\of{(L^p \cdot W)} \ge \sigma_p$ for any subvariety $W$ of
dimension $p$ which contains $x_0$,

\item $\sigma_1 \ge 3, \sigma_2 \ge 3$, and $\sigma_3 > 3$.
\endroster
Then $\vert K_X + L \vert$ is free at $x_0$.
\endproclaim

\demo{Proof}
{\it Step 0}.
Let $t$ be a rational number such that
$t > \frac{3}{\root{3}\of{(L^3)}}$.
Since $\sigma_3 > 3$, we can take $t < 1$.
By Proposition 2.1, there exists an effective $\Bbb Q$-Cartier
divisor $D$ such that
$D \sim_{\Bbb Q} tL$ and $\text{ord}_{x_0}D = 3$.
Let $c \le 1$ be the {\it log canonical threshold} of $(X, D)$ at $x_0$:
$$
c = \sup \{t \in \Bbb Q ; K_X + tD \text{ is LC at }x_0 \}
$$
and let $W$ be the minimal element of $CLC(X, x_0, cD)$.
If $W = \{x_0\}$, then $\vert K_X + L \vert$ is free at $x_0$
by Propositions 2.3.

\vskip .5pc

{\it Step 1}.
We consider the case in which $W = C$ is a curve.
Since $t < 1$, we have $ct + (1 - c) < 1$.
Since $\sigma_1 \ge 3$,
there exists a rational number $t'$ with $ct + (1-c) < t' < 1$ and
an effective $\Bbb Q$-Cartier divisor $D'_C$ on $C$ such that
$D'_C \sim_{\Bbb Q} (t'-ct)L \vert _C$ and $\text{ord}_{x_0}D'_C = 3(1-c)$.
By the Serre vanishing theorem, there exists
an effective $\Bbb Q$-Cartier divisor $D'$ on $X$ such that
$D' \sim_{\Bbb Q} (t'-ct)L$ and $D' \vert _C = D'_C$.
In fact, if we take a sufficiently large and divisible integer $m$
such that $mD'_C$ is a Cartier divisor in $\vert m(t'-ct)L \vert _C \vert$ and
$H^1(X, \Cal I_C(m(t'-ct)L)) = 0$, then
there exists an extension $D'_m \in \vert m(t'-ct)L \vert$ of $mD'_C$,
so we set $D' = \frac 1m D'_m$.

Let $D'_1$ be a general effective $\Bbb Q$-Cartier divisor
on an affine neighborhood $U$ of $x_0$ in $X$ such that
$D'_1 \vert_{C \cap U} = D'_C \vert_{C \cap U}$
and $\text{ord}_{x_0}D'_1 = 3(1-c)$.
Then we have $\text{ord}_{x_0}(cD + D'_1) = 3$, hence
$\{x_0\} \in CLC(U, cD + D'_1)$.
Let
$$
c' = \sup \{t \in \Bbb Q; K_X + (cD + tD'_1)
\text{ is LC at }x_0 \}.
$$
By Theorem 1.10, we conclude that $(X, cD + c'D')$
is LC at $x_0$, and $CLC(X, x_0, cD + c'D')$ has an element which
is properly contained in $C$, i.e., $\{x_0\}$.

\vskip .5pc

{\it Step 2-1}.
We consider the case in which $W = S$ is a surface.
$S$ is smooth or has a rational double point at $x_0$.
Let $d:= \text{mult}_{x_0} S = 1$ or $2$.
We assume first that $d = 1$.
As in Step 1,
we take a rational number $t'$,
an effective $\Bbb Q$-Cartier divisor $D'$ on $X$
and a positive number $c'$ such that $ct + (1-c) < t' < 1$,
$D' \sim_{\Bbb Q} (t'-ct)L$, $\text{ord}_{x_0}D' \vert_S = 3(1-c)$,
$(X, cD + c'D')$ is LC at $x_0$, and that
the minimal element $W'$ of $CLC(X, x_0, cD + c'D')$
is properly contained in $S$.
Thus we have the theorem when $W' = \{x_0\}$.

We consider the case in which $W' = C$ is a curve.
Since $t, t'< 1$, we have
$ct + c'(t'-ct) + (1-c)(1-c') < 1$.  As in Step 1,
there exists a rational number $t''$ with $ct + c'(t'-ct) + (1-c)(1-c')
< t'' < 1$ and
an effective $\Bbb Q$-Cartier divisor $D''_C$ on $C$ such that
$D''_C \sim_{\Bbb Q} (t''-ct-c'(t'-ct))L \vert_C$
and $\text{ord}_{x_0}D''_C = 3(1-c)(1-c')$.
Let $D'' \sim_{\Bbb Q} (t''-ct-c'(t'-ct))L$ be its extension to $X$ as before.
Let us take a general effective $\Bbb Q$-Cartier divisor $D''_1$
on an affine neighborhood $U$ of $x_0$ in $X$ such that
$D''_1 \vert_{C \cap U} = D''_C \vert_{C \cap U}$
and $\text{ord}_{x_0}D''_1 \vert_S = 3(1-c)(1-c')$.
By Theorem 1.10, there exists $c''$ such that
$(U, cD + c'D' + c''D''_1)$
is LC at $x_0$ and there exists an element of
$CLC(U, x_0, cD + c'D' + c''D''_1)$ which is properly contained in $S$.
Moreover, since $D''_1$ is chosen to be general, we have $c'' > 0$ and
the minimal element should be properly contained in $C$.
By Theorem 1.10 again, $(X, cD + c'D' + c''D'')$ is LC at $x_0$
and there exists an element of $CLC(X, x_0, cD + c'D' + c''D'')$
which is properly contained in $C$.

\vskip .5pc

{\it Step 2-2}.
We assume that $d = 2$.
As in Step 2-1,
we take a rational number $t'$ with $ct + \sqrt 2(1-c) < t'$ and
an effective $\Bbb Q$-Cartier divisor $D'$ on $X$ with
$D' \sim_{\Bbb Q} (t'-ct)L$ and
$\text{ord}_{x_0}D' \vert_S = 3(1-c)$.
Here we need the factor $\sqrt 2$ because $S$ has multiplicity 2 at $x_0$.
We take $0 < c' \le 1$ such that
$(X, cD + c'D')$ is LC and
$CLC(X, x_0, cD + c'D')$ has an element which is properly contained in $S$.

We shall prove that we may assume $ct + \sqrt 2(1-c) < 1$.  Then we can take
$t' < 1$ as in Step 2-1, and the rest of the proof is the same.
For this purpose, we apply Theorem 2.2.
In the argument of Steps 0 through 2-1,
the number $t$ was chosen under the only condition that $t < 1$.
So we can take $t = 1 - \epsilon_1$,
where the $\epsilon_n$ for $n = 1,2, \ldots$
will stand for very small positive rational numbers.
Then $k = \frac 3{1 - \epsilon_1} = 3 + \epsilon_2$ and $e = \frac1{3c}$.
This means the following: for any effective $D \sim_{\Bbb Q} tL$, if
$\text{ord}_{x_0} D \ge 3$, then $cD \ge S$.
We look for $k' = 6$ so that $D \sim_{\Bbb Q} tL$
with $t = \frac 12$ and $\text{ord}_{x_0} D \ge 3$.
The equation for $k'$ becomes
$$
\lambda^3 + (\frac{1 - 2e - \lambda}{1 - \lambda})^2
\{(2 + \frac{2\lambda e}{1 - 2e - \lambda})^3
- (\lambda + \frac{2\lambda e}{1 - 2e - \lambda})^3\} < 1.
$$
We have $\lambda + 2e \le 1$, $\frac 13 \le e$, and
in particular, $0 \le \lambda \le \frac 13$.
If we put $\frac {2e}{1-\lambda} = \alpha$, then
$\frac 2{3(1-\lambda)} \le \alpha \le 1$, and
$$
(- \lambda^3 + 6\lambda^2 - 12\lambda + 8)\alpha^2
+ (- \lambda^3 + 12\lambda - 16)\alpha + 8 < 1.
$$
For a fixed $\lambda$, since $- \lambda^3 + 6\lambda^2 - 12\lambda + 8 > 0$,
the left hand side attains the maximum at
$\alpha = \frac 2{3(1-\lambda)}$ or $= 1$, and the values
are given by
$\frac 1{(1-\lambda)^2}(\frac 23 \lambda^4 - \frac {10}9 \lambda^3
+ \frac 83 \lambda^2 - \frac 83 \lambda + \frac 89)$ or
$- 2\lambda^3 + 6\lambda^2$, respectively.
Since both numbers are smaller than $1$, we obtain a member of
$\vert mL \vert$ with $k' = 6$ for some $m$.
Then we choose a new $t$ as $t = \frac 12$.
We repeat the preceding argument from Step 0.
If we arrive at Step 2-2 again, then
we have $\frac23 \le c \le 1$ and
$ct + \sqrt{2}(1-c) < 1$.
\hfill $\square$
\enddemo

\proclaim{Corollary 3.2} ([EL1], [F]).
Let $X$ be a smooth projective variety of dimension $3$,
and $H$ an ample divisor.
Then $\vert K_X + mH \vert$ is free if $m \ge 4$.
Moreover, if $(H^3) \ge 2$, then $\vert K_X + 3H \vert$ is also free.
\hfill $\square$
\endproclaim

\head 4. Smooth 4-fold
\endhead

\proclaim{Theorem 4.1}
Let $X$ be a normal projective variety of dimension $4$,
$L$ an ample Cartier divisor, and $x_0 \in X$ a smooth point.
Assume that there are positive numbers $\sigma_p$ for $p = 1,2,3,4$
which satisfy the following conditions:

(1) $\root{p}\of{(L^p \cdot W)} \ge \sigma_p$ for any subvariety $W$ of
dimension $p$ which contains $x_0$,

(2) $\sigma_p \ge 5$ for all $p$.

Then $\vert K_X + L \vert$ is free at $x_0$.
\endproclaim

\demo{Proof}
{\it Steps 0, 1 and 2-1}.
Let $t$ be a rational number such that
$t > \frac{4}{\root{4}\of{(L^4)}}$.
Since $\sigma_4 > 4$, we can take $t < 1$.
By Proposition 2.1,
there exists an effective $\Bbb Q$-Cartier divisor $D$ such that
$D \sim_{\Bbb Q} tL$ and $\text{ord}_{x_0}D = 4$.
Let $c \le 1$ be the canonical threshold of $(X, D)$ at $x_0$, and
$W$ the minimal element of $CLC(X, x_0, cD)$.
If $W$ is a point, then $\vert K_X + L \vert$ is free at $x_0$
by Propositions 2.3.
Since $\sigma_p \ge 4$ for $p = 1,2$,
the cases in which $W$ is a curve or a smooth surface can be treated
similarly as in Steps 1 and 2-1 of the proof of Theorem 3.1.

\vskip .5pc

{\it Step 2-2}.
We consider the case in which $W = S$ is a surface.
Since $S$ has a rational singularity at $x_0$ and its embedding dimension
is 3 or 4, we have
$d := \text{mult}_{x_0} S = 2$ or $3$ ([A]).
We can take $t < \frac 45 + \epsilon$ for $0 < \epsilon \ll 1$, so
$\frac {3t}4 + \frac {\sqrt 3}5 < 1$.
Therefore, if $c \ge \frac 34$, then
$$
ct + \frac {4\sqrt d(1-c)}{\sigma_2} < 1  \tag a
$$
In this case, we can
take a rational number $t'$ and an effective $\Bbb Q$-Cartier divisor $D'$ on
$X$ such that $ct + \frac {4\sqrt d(1-c)}{\sigma_2} < t' < 1$,
$D' \sim_{\Bbb Q} (t'-ct)L$ and $\text{ord}_{x_0}D' \vert_S = 4(1-c)$,
and proceed as in Step 2-1.
On the other hand, if $c \le \frac 34$, then
$$
ct + \frac {\sqrt d}{\sigma_2} < 1  \tag b
$$
We take $t'$ and $D'$ with
$D' \sim_{\Bbb Q} (t'-ct)L$ and $\text{ord}_{x_0}D' \vert_S = 1$,
and use Theorem 1.9 in order to obtain a smaller
center of log canonical singularities.

Let $c'$ and $W'$ as before.
We consider the case in which $W' = C$ is a curve.
We have
$$
\text{ord}_{x_0}(cD + c'D') \vert_S \ge
\cases 4c + 4c'(1-c) = 4 - 4(1-c)(1-c') &\text{ in the case (a)} \\
4c + c' = 4 - 4(1-c-\frac{c'}4) &\text{ in the case (b)} \endcases
$$
In the case (a), since $t,t' < 1$, we have
$ct + c'(t'-ct) + (1-c)(1-c') < 1$.
In the case (b), we have
$$
ct + c'(t'-ct) + \frac 45 (1-c-\frac{c'}4)
< (\frac45 + \epsilon)c + \frac 45(1-c) +
c'(\frac {\sqrt 3}5 - \frac15)
\le \frac{3+\sqrt3}5 + \epsilon < 1
$$
The rest is the same as before.

\vskip .5pc

{\it Step 3}.
We consider the case in which $W = V$ is a 3-fold.
Let $d:= \text{mult}_{x_0} V$.
By Proposition 1.7, $V$ has only a canonical singularity at $x_0$, and
we have $d = 1, 2$ or $3$.
Since $t < \frac 45 + \epsilon$ and $c \ge \frac d4$, we have
$$
ct + \frac {4\root{3}\of{d}(1-c)}{\sigma_3}
< \frac {d + (4-d)\root{3}\of{d}}5 + \epsilon
< 1
$$
Therefore, there exists $t' < 1$ and $D'$ as before, and we obtain
$c'$ and $W'$.
If $W' = C$ is a curve, then
$$
ct + \frac {4\root{3}\of{d}(1-c)c'}{\sigma_3}
+ \frac {4(1-c)(1-c')}{\sigma_1} < 1
$$
and we obtain our assertion as before.

We assume that $W' = S$ is a surface.
Let $d':= \text{mult}_{x_0} S = 1, 2$ or $3$.
We have to prove that one of the followings hold:
$$
\align
&ct + \frac {4\root{3}\of{d}(1-c)c'}{\sigma_3}
+ \frac {4\sqrt {d'}(1-c)(1-c')}{\sigma_2} < 1 \tag a \\
&ct + \frac {4\root{3}\of{d}(1-c)c'}{\sigma_3}
+ \frac {\sqrt {d'}}{\sigma_2} < 1  \tag b
\endalign
$$
If $\root{3}\of{d} \ge \sqrt{d'}$, then (a) holds.  Otherwise,
we have the following cases:

\roster
\item $d = 1$, $d' = 2$ or $3$,

\item $d = 2$, $d' = 2$ or $3$,

\item $d = 3$, $d' = 3$.
\endroster
In the case (1), since the embedding dimension of $S$ at $x_0$ is 3,
we have $d' = 2$.  Then we have $4c \ge 1$ and $4(1-c)c' \ge 2$, hence
$$
\text{l.h.s. of (a)} <
\frac {1+2+\sqrt2}5 + \epsilon < 1
$$
In the case (2), since $4c \ge 2$, if $4(1-c)c' \ge 1$, then
$$
\text{l.h.s. of (a)} <
\frac {2+\root{3}\of{2}+\sqrt3}5 + \epsilon
= \frac{4.9919 \cdots}5 + \epsilon < 1
$$
Otherwise, (b) holds.
Finally, in the case (3),
$$
\text{l.h.s. of (a)} <
\frac {3+\root{3}\of{3}c'+\sqrt3(1-c')}5 + \epsilon <
\frac {3+\sqrt3}5 + \epsilon < 1
$$
In any case, we obtain $t'' < 1$, $D''$,
$c''$ and $W''$ as before.

When $W''$ is a curve, we have still $t''' < 1$.
In fact, we have
$$
ct + \frac {4\root{3}\of{d}(1-c)c'}{\sigma_3}
+ \frac {4\sqrt {d'}(1-c)(1-c')c''}{\sigma_2}
+ \frac {4(1-c)(1-c')(1-c'')}{\sigma_1}< 1
$$
in the case (a), and
$$
ct + \frac {4\root{3}\of{d}(1-c)c'}{\sigma_3}
+ \frac {\sqrt {d'}c''}{\sigma_2}
+ \frac {4(1-c)(1-c') - c''}{\sigma_1} < 1
$$
in the case (b).
The rest of the proof is similar to the previous steps.
\hfill $\square$
\enddemo

\proclaim{Corollary 4.2}
Let $X$ be a smooth projective variety of dimension $4$,
and $H$ an ample divisor.
Then $\vert K_X + mH \vert$ is free if $m \ge 5$.
\hfill $\square$
\endproclaim


\Refs
\widestnumber\key{KMM}

\ref\key A
\by M. Artin
\paper On isolated rational singularities of surfaces
\jour Amer. J. Math. \vol 88 \pages 129--136
\yr 1966
\endref

\ref\key AS
\by U. Angehrn and Y.-T. Siu
\paper Effective freeness and separation of points for adjoint bundles
\jour Invent. Math
\paperinfo to appear
\endref

\ref\key EL1
\by L. Ein and R. Lazarsfeld
\paper Global generation of pluricanonical and adjoint linear series
on smooth projective threefolds
\jour J. Amer. Math. Soc. \vol 6 \pages 875--903
\yr 1993
\endref

\ref\key EL2
\bysame
\paperinfo manuscript
\endref

\ref\key F
\by T. Fujita
\paper Remarks on Ein-Lazarsfeld criterion of spannedness of adjoint bundles of
polarized threefold
\paperinfo preprint
\endref


\ref\key K1
\by Y. Kawamata
\paper A generalization of Kodaira-Ramanujam's vanishing theorem
\jour Math. Ann. \vol 261 \pages 43--46
\yr 1982
\endref

\ref\key K2
\bysame
\paper On the finiteness of generators of a pluri-canonical ring for a 3-fold
of general type
\jour Amer. J. Math. \vol 106 \pages 1503--1512
\yr 1984
\endref

\ref\key K3
\bysame
\paper Subadjunction of log canonical divisors for a subvariety of
codimension 2
\paperinfo \break preprint
\endref

\ref\key KMM
\by Y. Kawamata, K. Matsuda and K. Matsuki
\paper Introduction to the minimal model problem
\jour Adv. St. Pure Math. \vol 10 \pages 283--360
\yr 1987
\endref

\ref\key Ko
\by J. Koll\'ar et.al.
\paper Flips and abundance for algebraic threefolds
\jour Ast\'erisque \vol 211
\yr 1992
\endref

\ref\key Sh1
\by V. V. Shokurov
\paper The non-vanishing theorem
\jour Math. USSR-Izv. \vol 26 \pages 591--604
\yr 1986
\endref

\ref\key Sh2
\bysame
\paper 3-fold log flips
\jour Math. USSR-Izv. \vol 56 \pages 105--203
\yr 1992
\endref

\ref\key Sm
\by K. Smith
\paper Fujita's freeness conjecture in terms of local cohomology
\paperinfo preprint
\yr 1995
\endref

\ref\key T
\by H. Tsuji
\paper Global generation of adjoint bundles
\paperinfo preprint
\yr 1994
\endref

\ref\key V
\by E. Viehweg
\paper Vanishing theorems
\jour J. reine angew. Math. \vol 335 \pages 1--8
\yr 1982
\endref

\endRefs
\enddocument